\documentclass{amsproc}

\usepackage{latexsym}
\usepackage{amsmath}
\usepackage{amsfonts}
\usepackage{amssymb}
\usepackage{amscd}

\numberwithin{equation}{section}

\newtheorem{proposition}{Proposition}[section]
\newtheorem{lemma}[proposition]{Lemma}
\newtheorem{theorem}[proposition]{Theorem}
\newtheorem{corollary}[proposition]{Corollary}

\newtheorem{definition}[proposition]{Definition}

\newtheorem{remark}[proposition]{Remark}

\def\itheorem#1#2{\newtheorem{#1}[proposition]{#2}}

\renewcommand{\labelenumi}{(\roman{enumi})}

\itheorem{Remark}{Remark}

\def\Hom{\mathop{\rm Hom}\nolimits}
\def\Ext{\mathop{\rm Ext}\nolimits}

\def\cf{\mathop{\rm cf}\nolimits}
\def\im{\mathop{\rm Im}\nolimits}
\def\dom{\mathop{\rm Dom}\nolimits}
\def\rg{\mathop{\rm rg}\nolimits}
\def\pr{\mathop{\rm pr}\nolimits}

\def\rk{\mathop{\rm rk}}

\def\cl{\mathop{\rm cl}}

\def\P{{\mathcal{P}}}

\def\v{\mathcal{V}}

\def\fk{\mathop{\rm fr-rk}\nolimits}

\def\Q{{\mathbb{Q}}}
\def\G{{\mathbb{G}}}
\def\Z{{\mathbb{Z}}}
\def\N{{\mathbb{N}}}

\begin{document}

\title{On the $p$-rank of $\Ext_{\Z}(G,\Z)$ in certain models of $ZFC$}

\author{Saharon Shelah}

\address{Department of Mathematics,
The Hebrew University of Jerusalem, Israel, and Rutgers
University, New Brunswick, NJ U.S.A.}

\email{Shelah@math.huji.ac.il}
\thanks{Number 874 in Shelah's list of publications. The first author
  was supported by project No. I-706-54.6/2001 of the {\em
    German-Israeli Foundation for Scientific Research \&
    Development}.\\ The second author was supported by a grant from
  the German Research Foundation DFG}

\author{Lutz Str\"ungmann}

\address{Department of Mathematics,
University of Duisburg-Essen, 45117 Essen, Germany}

\email{lutz.struengmann@uni-essen.de}
\address{{\it Current address}: Department of Mathematics, University of Hawaii, 2565 McCarthy Mall, Honolulu, HI 96822-2273, USA}

\email{lutz@math.hawaii.edu}

\begin{abstract}
We show that if the existence of a supercompact cardinal is
consistent with $ZFC$, then it is consistent with $ZFC$ that the
$p$-rank of $\Ext_{\Z}(G,\Z)$ is as large as possible for every
prime $p$ and any torsion-free abelian group $G$. Moreover, given an
uncountable strong limit cardinal $\mu$ of countable cofinality and
a partition of $\Pi$ (the set of primes) into two disjoint subsets
$\Pi_0$ and $\Pi_1$, we show that in some model which is very close
to $ZFC$ there is an almost-free abelian group $G$ of size
$2^{\mu}=\mu^+$ such that the $p$-rank of $\Ext_{\Z}(G,\Z)$ equals
$2^{\mu}=\mu^+$ for every $p \in \Pi_0$ and $0$ otherwise, i.e. for
$p \in \Pi_1$.
\end{abstract}

\maketitle \setcounter{section}{0}

\section{Introduction}
In 1977 the first named author solved the well-known Whitehead problem by showing that it is
undecidable in ordinary set-theory $ZFC$ wether or not every abelian group
$G$ satisfying $\Ext_{\Z}(G,\Z)=\{0\}$ has to be free (see \cite{Sh1}, \cite{Sh2}).
However, this did not clarify the structure of $\Ext_{\Z}(G,\Z)$ for
torsion-free abelian groups - a problem which has received much attention since then.
Easy arguments show that $\Ext_{\Z}(G,\Z)$ is always a divisible group for
every torsion-free group $G$. Hence it is of the form
\[ \Ext_{\Z}(G,\Z)= \bigoplus\limits_{p \in
\Pi}\Z(p^{\infty})^{(\nu_p)} \oplus \Q^{(\nu_0)} \] for some
cardinals $\nu_p,\nu_0$ ($p \in \Pi)$ which are uniquely determined.
The obvious question that arises is which sequences $(\nu_0,
\nu_p : p \in \Pi)$ of cardinals can appear as the cardinal invariants of
$\Ext_{\Z}(G,\Z)$ for some (which) torsion-free abelian group? Obviously,
the trivial sequence consisting of zero entries only can be realized
by any free abelian group. However, the solution of the Whitehead
problem shows that it is undecidable in $ZFC$ if these are the only
ones. There are a few results about possible sequences $(\nu_0,
\nu_p : p \in \Pi)$ provable in $ZFC$. On the other hand, assuming
G\"odel's constructible universe $(V=L)$ plus there is no weakly compact cardinal
a complete characterization of the cardinal
invariants of $\Ext_{\Z}(G,\Z)$ for torsion-free abelian groups $G$
has recently been completed by the authors (see \cite{EkHu},
\cite{EkSh}, \cite{GS}, \cite{GS2}, \cite{HHS}, \cite{MRS}, \cite{SS1},
\cite{SS2}, \cite{ShSt} and \cite{Sh3} for references). In fact, it turned out that
almost all divisible groups $D$ may be realized as $\Ext_{\Z}(G,\Z)$ for some torsion-free
abelian group $G$ of almost any given size.\\
In this paper we shall take the opposite point of view. It is a theorem
of $ZFC$ that every sequence $(\nu_0,\nu_p : p \in \Pi)$ of cardinals such that
$\nu_0=2^{\lambda_0}$ for some infinite $\lambda_0$ and $\nu_p\leq \nu_0$ is either finite
or of the form $2^{\lambda_p}$ for some infinite $\lambda_p$ can arise as the
cardinal invariants of $\Ext_{\Z}(G,\Z)$ for some torsion-free $G$. The first purpose
of this paper is to show that this result is as best as possible by constructing a model
of $ZFC$ in which the only realizable sequences are of this kind. We shall
assume therefore the consistency of the existence of a supercompact cardinal (see \cite{MeSh}). This is a strong
additional set-theoretic assumption which makes the model we are working in be far from $ZFC$.\\
On the other hand, we will also work in models very close to $ZFC$ assuming only the existence
of certain ladder systems on successor of strong limit cardinals of cofinality $\aleph_0$.
Although this model is close to $ZFC$ it allows us to construct almost-free torsion-free abelian groups
$G$ such that for instance $\Ext_{\Z}(G,\Z)$ is torsion-free, i.e. $G$ is coseparable. Also this can be considered
as a result at the borderline of what is provable in models close to $ZFC$ since
the existence of non-free coseparable groups is independent of $ZFC$ (see \cite[Chapter XII]{EkMe} and \cite{MeSh}).\\

Our notation is standard and we write maps from the right. All groups under consideration are abelian and written additively.
We shall abbreviate $\Ext_{\Z}(-,-)$ by $\Ext(-,-)$ and $\Pi$ will
denote the set of all primes. A Whitehead group is a torsion-free
group $G$ such that $\Ext_{\Z}(G,\Z)=0$. If $H$ is a pure subgroup of the
abelian group $G$, then we shall write $H \subseteq_* G$. We shall assume
sufficient knowledge about forcing, large cardinals and prediction principles
like weak diamond etc. as for example in \cite{EkMe}, \cite{Ku} or \cite{Sh4}.
Also reasonable knowledge about abelian groups as for instance in \cite{Fu} is
assumed. However, the authors have tried to make the paper as accessible as possible
to both algebraists and set theorists.

\section{The structure of $\Ext(G,\Z)$}
In this section we recall the basic results on the structure of
$\Ext(G,\Z)$ for torsion-free groups $G$. Therefore let $G$ be a
torsion-free abelian group. It is easy to see that $\Ext(G,\Z)$ is
divisible, hence it is of the form \[ \Ext(G,\Z)=
\bigoplus\limits_{p \in \Pi}\Z(p^{\infty})^{(\nu_p)} \oplus
\Q^{(\nu_0)} \] for certain cardinals $\nu_p,\nu_0$ ($p \in \Pi)$.
Since the cardinals $\nu_p$ ($p \in \Pi)$ and $\nu_0$ completely
determine the structure of $\Ext(G,\Z)$ we introduce the following
terminology. Let $\Ext_p(G,\Z)$ be the $p$-torsion part of
$\Ext(G,\Z)$ for $p \in \Pi$. We denote by $r^e_0(G)$ the {\it
torsion-free rank} $\nu_0$ of $\Ext(G,Z)$ which is the dimension of
$\Q \otimes \Ext(G,\Z)$ and by $r_p^e(G)$ the $p$-$rank$ $\nu_p$ of
$\Ext(G,\Z)$ which is the dimension of $\Ext(G,\Z)[p]$ as a vector
space over $\Z/p\Z$ for any prime number $p \in \Pi$. There are only
a few results provable in $ZFC$ when $G$ is uncountable, but
assuming additional set-theoretic assumptions a better understanding
of the structure of $\Ext(G,\Z)$ is obtained. For instance, in
G\"odel's universe and assuming that there is no weakly compact
cardinal a complete characterization is known. The aim of this paper
is to go to the borderline of the characterization. On the one hand
we shall show that one can make the $p$-ranks of $\Ext(G,\Z)$ as
large as possible for every torsion-free abelian group $G$ by
working in a model of $ZFC$ which assumes strong additional axioms
(the existence of large cardinals). On the other hand we shall work
in a model which is very close to $ZFC$ but still allows
to construct uncountable torsion-free groups $G$ such that $\Ext(G,\Z)$ is torsion-free.\\
We first justify our restriction to torsion-free $G$. Let $A$ be any
abelian group and $t(A)$ its torsion subgroup. Then
$\Hom(t(A),\Z)=0$ and hence we obtain the short exact sequence
\[0 \rightarrow \Ext(A/t(A),\Z) \rightarrow \Ext(A,\Z)
\rightarrow \Ext(t(A),\Z) \rightarrow 0 \] which must split since
$\Ext(A/t(A),\Z)$ is divisible. Thus \[ \Ext(A,\Z) \cong
\Ext(A/t(A),\Z) \oplus \Ext(t(A),\Z). \] Since the structure of
$\Ext(t(A),\Z) \cong \prod_{p \in \Pi} \Hom(A,\Z(p^{\infty}))$ is
well-known in $ZFC$ (see \cite{Fu}) it is reasonable to assume that
$A$ is torsion-free and, of course, non-free. Using Pontryagin's
theorem one proves

\begin{lemma}
\label{countabletf} Suppose $G$ is a countable torsion-free group
which is not free. Then $r_0^e(G)=2^{\aleph_0}$.
\end{lemma}

\proof See \cite[Theorem XII 4.1]{EkMe}. \qed\\

Similarly, we have for the $p$-ranks of $G$ the following

\begin{lemma}
\label{countablet} If $G$ is a countable torsion-free group, then
for any prime $p$, either $r_p^e(G)$ is finite or $2^{\aleph_0}$.
\end{lemma}

\proof See \cite[Theorem XII 4.7]{EkMe}. \qed\\

This clarifies the structure of $\Ext(G,\Z)$ for countable
torsion-free groups $G$ in $ZFC$. We now turn our attention to
uncountable groups. There is a useful characterization of $r_p^e(G)$
using the exact sequence
\[ 0 \rightarrow \Z \overset{p}\rightarrow \Z \rightarrow \Z/pZ \rightarrow
0.
\]
The induced sequence
\[ \Hom(G,\Z) \overset{\varphi^p}{\rightarrow} \Hom(G,\Z/p\Z)
\rightarrow \Ext(G,\Z) \overset{p_*}{\rightarrow} \Ext(G,\Z) \]
shows that the dimension of
\[ \Hom(G,\Z/p\Z)/\Hom(G,\Z)\varphi^p \]
as a vector space over $\Z/p\Z$ is exactly $r_p^e(G)$.\\

The following result due to Hiller, Huber and Shelah deals with the
case when $\Hom(G,\Z)=0$.

\begin{lemma}
\label{existence} For any cardinal $\nu_0$ of the form
$\nu_0=2^{\mu_0}$ for some infinite $\mu_0$ and any sequence of
cardinals $(\nu_p : p \in \Pi)$ less than or equal to $\nu_0$ such
that each $\nu_p$ is either finite or of the form $2^{\mu_p}$ for
some infinite $\mu_p$ there is a torsion-free group $G$ of
cardinality $\mu_0$ such that $\Hom(G,\Z)=0$ and $r_0^e(G)=\nu_0$,
$r_p^e(G)=\nu_p$ for all primes $p \in \Pi$.
\end{lemma}

\proof See \cite[Theorem 3(b)]{HHS}. \qed\\

Together with the following lemma we have reached the borderline of
what is provable in $ZFC$.

\begin{lemma}
\label{dualt} If $G$ is torsion-free such that $\Hom(G,\Z)=0$, then
for all primes $p$, $r_p^e(G)$ is either finite or of the form
$2^{\mu_p}$ for some infinite $\mu_p \leq |G|$. \end{lemma}

\proof See \cite[Lemma XII 5.2]{EkMe}. \qed\\

Assuming G\"odel's axiom of constructibility one even knows a
complete characterization in the case when $\Hom(G,\Z)=0$.

\begin{lemma}[V=L]
\label{uncountabletf} Suppose $G$ is a torsion-free non-free group
and let $B$ be a subgroup of $A$ of minimum cardinality $\nu$ such
that $A/B$ is free. Then $r_0^e(G)=2^{\nu}$. In particular,
$r_0^e(G)$is uncountable and $r_0^e(G)=2^{|G|}$ if $\Hom(G,\Z)=0$.
\end{lemma}

\proof See \cite[Theorem XII 4.4, Corollary XII 4.5]{EkMe}. \qed\\

Note that the above lemma is not true in $ZFC$ since for any
countable divisible group $D$ it is consistent that there exists an
uncountable torsion-free group $G$ with $\Ext(G,\Z) \cong D$,
hence $r_0^e(G)=1$ is possible taking $D=\Q$ (see \cite{Sh3}).\\

The following result is a collection of theorems due to Grossberg,
Mekler, Roslanowski, Sageev and the authors. It shows that under the
assumption of $(V=L)$ almost all possibilities for $r_p^e(G)$ can
appear if the group is not of weakly compact cardinality or singular
cardinality of cofinality $\aleph_0$.

\begin{lemma}[V=L]
Let $\nu$ be an uncountable cardinal and suppose that $(\nu_p : p\in
\Pi)$ is a sequence of cardinals such that for each $p$, $0 \leq
\nu_p \leq 2^{\nu}$. Moreover, let $H$ be a torsion-free group of
cardinality $\nu$. Then the following hold.

\begin{enumerate}
\item If $\nu$ is regular and less than the first weakly compact cardinal, then there is an almost-free group $G$ of cardinality $\nu$
such that $r_0^e(G)=2^{\nu}$ and for all primes $p$, $r_p^e(G)=\nu_p$;
\item If $\nu$ is a singular strong limit
cardinal of cofinality $\omega$, then there is no torsion-free group
$G$ of cardinality $\nu$ such that $r_p^e(G)=\nu$ for any prime $p$;
\item If $\nu$ is weakly compact and $r_p^e(H) \geq \nu$ for some prime $p$, then
$r_p^e(H)=2^{\nu}$; \item If $\nu$ is singular less than the first
weakly compact cardinal and of cofinality $\cf(\nu) > \aleph_0$,
then there is a torsion-free group $G$ of cardinality $\nu$ such
that $r_0^e(G)=2^{\nu}$ and for all primes $p$, $r_p^e(G)=\nu_p$.
\end{enumerate}
\end{lemma}

\proof For (i) see \cite[Theorem 3.7]{MRS}, for (ii) we refer to \cite[Theorem 1.0]{GS},
for (iii) see \cite[Main Theorem]{SS1} and (iv) is contained in \cite{ShSt}. \qed\\

The above results show that under the assumption of $(V=L)$ and the
non-existence of weakly compact cardinals, the structure of
$\Ext(G,\Z)$ for torsion-free groups $G$ of cardinality $\nu$ is
clarified for all cardinals $\nu$ and almost all sequences
$(\nu_0,\nu_p : p \in \Pi)$ can be realized as the cardinal
invariants of some torsion-free abelian group in almost every
cardinality. However, if we weaken the set-theoretic assumptions to
$GCH$ (the generalized continuum hypothesis), then even more is
possible which was excluded by $(V=L)$ before (see Lemma
\ref{uncountabletf}).

\begin{lemma}The following hold.
\begin{enumerate}
\item Assume $GCH$. For any torsion-free group $A$ of uncountable
cardinality $\nu$, if $\Hom(A,\Z)=0$ and $r_0^e(A) < 2^{\nu}$, then
for each prime $p$, $r_p^e(A)=2^{\nu}$;
\item It is consisitent with $ZFC$ and $GCH$ that for any cardinal $\rho
\leq \aleph_1$, there is a torsion-free group $G_{\rho}$ such that
$\Hom(G_{\rho},\Z)=0$, $r_0^e(G_{\rho})=\rho$ and for all primes
$p$, $r_p^e(G_{\rho})=2^{\aleph_1}$.
\end{enumerate}
\end{lemma}

\proof See \cite[Theorem XII 5.3]{EkMe} and \cite[Theorem XII
5.49]{EkMe}.\qed\\

It is our aim in the next section to show that this rich structure
of $\Ext(G,\Z)$ ($G$ torsion-free) which exists in $(V=L)$ does not
appear in other models of $ZFC$. As a motivation we state two
results from \cite{MeSh} which show that using Cohen forcing we may
enlarge the $p$-rank of $\Ext(G,\Z)$ for torsion-free groups $G$.

\begin{lemma}
\label{meklerold} Suppose $G$ is contained in the $p$-adic
completion of a free group $F$ and $|G| > |F|$. Then, if $\lambda
\geq |F|$ and $\lambda$ Cohen reals are added to the universe,
$|\Ext_p(G,\Z)| \geq \lambda$. In particular, adding $2^{\aleph_0}$
Cohen reals to the universe implies that for every torsion-free
reduced non-free abelian group $G$ of cardinality less than the
continuum, there is a prime $p$ such that $r_p^e(G) > 0$.
\end{lemma}

\proof See \cite[Theorem 8]{MeSh}. \qed\\

Assuming the consistency of large cardinals we even get more. Recall
that a cardinal $\kappa$ is {\it compact} if it is uncountable
regular and satisfies the condition that for every set $S$, every
$\kappa$-complete filter on $S$ can be extended to a
$\kappa$-complete ultrafilter on $S$. This is equivalent to saying
that for any set $A$ such that $|A| \geq \kappa$, there exists a
fine measure on $P_{\kappa}(A)$ (the set of all subsets of $A$ of
size less than or equal to $\kappa$). If we require the measure to
satisfy a normality condition, then we get a stronger notion. A fine
measure $U$ on $P_{\kappa}(A)$ is called {\it normal} if whenever
$f:P_{\kappa}(A) \rightarrow A$ is such that $f(P) \in P$ for almost
all $P \in P_{\kappa}(A)$, then $f$ is constant on a set in $U$. A
cardinal $\kappa$ is called {\it supercompact} if for every set $A$
such that $|A| \geq \kappa$, there exists a normal measure on
$P_{\kappa}(A)$ (see \cite[Chapter II.2]{EkMe} or \cite[Chapter 6,
33. Compact cardinals]{Je} for further details on supercompact
cardinals).

\begin{lemma}
Suppose that it is consistent that a supercompact cardinal exists. Then it is consistent with either $2^{\aleph_0}=2^{\aleph_1}$ or $2^{\aleph_0} < 2^{\aleph_1}$ that for any group $G$ either $\Ext(G,\Z)$ is finite or $r_0^e(G) \geq 2^{\aleph_0}$.
\end{lemma}

\proof See \cite[Theorem 11]{MeSh}. \qed\\

\section{The free (p-)rank}
In this section we introduce the {\it free (p-)rank} of a torsion-free group $G$
($p$ a prime) which will induce upper bounds for the cardinal
invariants of $\Ext(G,\Z)$.

\begin{definition}
For a prime $p \in \Pi$ let $K_p$ be the class of all torsion-free
groups $G$ such that $G/p^{\omega}G$ is free. Moreover, let $K_0$ be
the class of all free groups.
\end{definition}

Note, that for $G \in K_p$ ($p$ a prime) we have $G=p^{\omega}G
\oplus F$ for some free group $F$ and hence $\Ext(G,\Z)[p]=0$ since
$p^{\omega}G$ is $p$-divisible. Note that $p^{\omega}G$ is a pure
subgroup of $G$. Thus $r_p^e(G)=0$ for $G \in K_p$ and any prime
$p$. Clearly, also $r_0^e(G)=0$ for all $G \in K_0$.

\begin{definition}
Let $G$ be a torsion-free group. We call
\[\fk_0(G)=\min\{\rk(H) : H \subseteq_* G \textit{ such that } G/H \in K_0 \} \] the {\rm free rank of $G$} and similarly we call
\[\fk_p(G)=\min\{\fk(H) : H \subseteq_* G/p^{\omega}G \textit{ such that } (G/p^{\omega}G)/H \in K_p \} \]
the {\rm free $p$-rank of $G$} for any prime $p \in \Pi$.
\end{definition}

We have a first easy lemma.

\begin{lemma}
\label{freerank}
Let $G$ be a torsion-free group and $p \in \Pi$ a prime. Then the following
hold.
\begin{enumerate}
\item $r_p^e(G)=r_p^e(G/p^{\omega}G)$;
\item If $H$ is a pure subgroup of $G$, then $r_p^e(H) \leq
r_p^e(G)$ and $r_0^e(H) \leq r_0^e(G)$;
\item $\fk_0(G)\geq \fk_0(G/p^{\omega}G)$;
\item $\fk_p(G)=\fk_p(G/p^{\omega}G)$;
\item $\fk_p(G) \leq \fk_0(G)$.
\end{enumerate}
\end{lemma}

\proof We first show (i) and let $p$ be a prime. Since $p^{\omega}G$
is pure in $G$ we have that $p^{\omega}G$ is $p$-divisible, hence
\[ 0 \rightarrow p^{\omega}G \rightarrow G \rightarrow G/p^{\omega}G
\rightarrow 0 \] induces the exact sequence
\[ 0 \rightarrow \Hom(G/p^{\omega}G,\Z/p\Z) \rightarrow
\Hom(G,\Z/p\Z) \rightarrow \Hom(p^{\omega}G,\Z/p\Z)=0 \] the latter
being trivial because $p^{\omega}G$ is $p$-divisible. Thus we have
\[ \Hom(G,\Z/p\Z) \cong \Hom(G/p^{\omega}G, \Z/p\Z) \] and it follows
easily that
\[ \Hom(G,\Z/p\Z)/\Hom(G,\Z)\varphi^{p}
\cong \Hom(G/p^{\omega}G,\Z/p\Z)/\Hom(G/p^{\omega}G,\Z)\varphi^{p}.
\]
Therefore, $r_p^e(G)=r_p^e(G/p^{\omega}G)$.\\
In order to show (ii) we consider the exact sequence
\[ 0 \rightarrow H \rightarrow G \rightarrow G/H \rightarrow 0 \]
which implies the exact sequence
\[ \cdots \rightarrow \Ext(G/H,\Z) \overset{\alpha}{\rightarrow} \Ext(G,\Z) \rightarrow
\Ext(H,\Z) \rightarrow 0. \] Since $G/H$ is torsion-free we conclude
that $\Ext(G/H,\Z)$ is divisible and hence $\im(\alpha)$ is
divisible. Thus $\Ext(G,\Z)= \Ext(H,\Z) \oplus \im(\alpha)$ and
therefore $r_0^e(H) \leq
r_0^e(G)$ and $r_p^e(H) \leq r_p^e(G)$ for every prime $p$.\\
 Claim (iii) is easily proved noting that,
whenever $G=H \oplus F$ for some free group $F$, then $p^{\omega}G
\subseteq H$ for every prime
$p$, hence $G/p^{\omega}G=H/p^{\omega}G \oplus F$.\\
To show (iv) note that $p^{\omega}(G/p^{\omega}G)=\{0\}$ and hence
$\fk_p(G)=\fk_p(G/p^{\omega}G)$ easily follows.\\
Finally, (v) follows from (iii), (iv) and the definition of
$\fk_p(G)$ and $\fk_0(G)$ since the class $K_0$ is contained in the
class $K_p$ for every prime $p$.
\qed\\

\begin{remark}
If $G$ is a torsion-free group and $p$ a prime, then Lemma \ref{freerank} (i) and (iv) imply that, regarding the free $p$-rank of $G$, we may assume
without loss of generality that $G$ is $p$-reduced. This is also justified by the fact that
\[ \fk_p(G)=\{ H \subseteq_* G \textit{ such that } G/H \in K_p \} \]
is easily proven.
\end{remark}

To simplify notations we let $\Pi_0=\Pi \cup \{0\}$ in the sequel.

\begin{lemma}
\label{lemmaabschaetz} Let $G$ be a torsion-free group. Then the following
hold.

\begin{enumerate}
\item $r_0^e(G) \leq 2^{\lambda}$ where $\lambda=\max\{\aleph_0,
\fk_0(G)\}$; In particular, $r_0^e(G) \leq 2^{\fk_0(G)}$ if
$\fk_0(G)$ is infinite;
\item $r_p^e(G) \leq p^{\fk_p(G)}$ for all $p \in \Pi$.
\end{enumerate}
\end{lemma}

\proof In order to prove (i) choose a subgroup $H \subseteq G$ such
that $\rk(H)=\fk_0(G)$ and $G/H \in K_0$. Hence $G=H \oplus F$ for
some free group $F$ and so $\Ext(G,\Z)=\Ext(H,\Z)$ which implies
that $r_0^e(G)=r_0^e(H)
\leq 2^{\lambda}$ where $\lambda=\max\{\aleph_0,\rk(H)\}=\max\{\aleph_0,\fk_0(G)\}$.\\
We now prove (ii). Let $p$ be a prime, then
$r_p^e(G)=r_p^e(G/p^{\omega}G)$ and $\fk_p(G)=\fk_p(G/p^{\omega}H)$
by Lemma \ref{freerank} (i) and (iv). Hence we may assume that
$p^{\omega}G=\{0\}$ without loss of generality. Let $H \subseteq G$
be such that $\fk_0(H)=\fk_p(G)$ and $G/H \in K_p$. Then $G/H=D
\oplus F$ for some free group $F$ and some $p$-divisible group $D$.
As in the proof of Lemma \ref{freerank} (i) it follows that
$r_p^e(G)=r_p^e(H)$. Now, we let $H=H' \oplus F'$ for some free
group $F'$ such that $\rk(H')=\fk_0(H)=\fk_p(G)$. Hence
$\Ext(H,\Z)=\Ext(H',\Z)$ and therefore
$r_p^e(G)=r_p^e(H)=r_p^e(H')$. Consequently, $r_p^e(G)=r_p^e(H')
\leq p^{\rk(H')}=p^{\fk_p(G)}$. \qed\\

Note, that for instance in $(V=L)$ for any torsion-free group $G$,
$2^{\fk_0(G)}$ is the actual value of $r_0^e(G)$ by Lemma
\ref{uncountabletf}. The following lemma justifies that, as far as
it concerns the free $p$-rank of a torsion-free group, one may also
assume without loss of generality that $\fk_p(G)=\rk(G)$ if $p \in
\Pi_0$.

\begin{lemma}
\label{reduction}
Let $G$ be a torsion-free group, $p \in \Pi_0$ and $H \subseteq_* G$ such that
\begin{enumerate}
\item $G/\left( H \oplus F \right) \in K_p$ for some free group $F$;
\item $\rk(H)=\fk_p(G)$.\end{enumerate}
Then $\fk_p(H)=\rk(H)$ and $r_p^e(G)=r_p^e(H)$.
\end{lemma}

\proof Let $G$, $H$ and $p$ be given. If $p=0$, then the claim is
trivially true. Hence assume that $p \in \Pi$ and that
$\rk(H)=\fk_p(G)$. Then there is a free group $F$ such that $H'=H
\oplus F$ is a pure subgroup of $G$ satisfying $G/H' \in K_p$. Thus
$\fk_p(G)=\rk(H)=\fk_0(H')$. Without loss of generality we may
assume that $G/H'$ is $p$-divisible by splitting of the free part.
By way of contradiction assume that $\fk_p(H) < \rk(H)$. Let $H_1
\subseteq_* H$ such that $H/H_1 \in K_p$ and $\fk_0(H_1)=\fk_p(H) <
\rk(H)$. Then there are a free group $F_1$ and a $p$-divisible group
$D$ such
\[ H/H_1 = D \oplus F_1. \]
Choose a pure subgroup $H_1 \subseteq_* H_2 \subseteq_* H$ such that
$H_2/H_1 \cong D$. Thus $H/H_2 \cong F_1$ and so $H \cong H_2 \oplus
F_1$ and without loss of generality $H=H_2 \oplus F_1$.
Consequently, $\rk(H_2)=\rk(H)$ since $\rk(H)=\fk_p(G)$. Let
$H_3=H_1 \oplus F_1 \oplus F$. Then
\[ \fk_0(H_3)=\fk_0(H_1) < \rk(H)=\fk_p(G). \]
Moreover,
\[ G/H' \cong \left( G/H_3 \right) / \left( H'/H_3 \right) \]
is $p$-divisible. Since also $H'/H_3 \cong H_2/H_1$ is $p$-divisible
and all groups under consideration are torsion-free we conclude that
$G/H_3$ is $p$-divisible. Hence $\fk_p(G) \leq \fk_0(H_3) <
\fk_p(G)$ - a contradiction. Finally, $r_p^e(G)=r_p^e(H)$ follows as
in the proof of Lemma \ref{lemmaabschaetz}. \qed\\

We now show how to calculate explicitly $\fk_p(G)$ for torsion-free
groups $G$ of finite rank and $p \in \Pi$ (note that $\fk_0(G)$ can
be easily calculated). Recall that a torsion-free group $G$ of
finite rank is {\it almost-free} if every subgroup $H$ of $G$ of
smaller rank than the rank of $G$ is free.

\begin{lemma}
Let $G$ be a non-free torsion-free group of finite rank $n$ and $p
\in \Pi$. Then we can calculate $\fk_p(G)$ as follows.

\begin{enumerate}
\item If $G$ is almost-free, then let $H \subseteq \Q$ be the outer type of $G$. Then
\begin{enumerate}
\item $\fk_p(G)=n$ if $H$ is not $p$-divisible;
\item $\fk_p(G)=0$ if $H$ is $p$-divisible.
\end{enumerate}
\item If $G$ is not almost-free, then choose a filtration $\{0\}=G_0
\subseteq_* G_1 \subseteq_* \cdots \subseteq_*G_m \subseteq_* G$
with $G_{k+1}/G_k$ almost-free. Then
\[ \fk_p(G)=\sum\limits_{k < m} \fk_p(G_{k+1}/G_k). \]
\end{enumerate}
\end{lemma}

\proof Left to the reader. \qed\\

In order to prove our main Theorem \ref{main1} of Section $4$ we
need a further result on the class $K_p$ for $p \in \Pi$.

\begin{lemma}
\label{forcingone}
Let $p$ be a prime and $G$ a torsion-free group of infinite rank. Then the following
hold.
\begin{enumerate}
\item If $G$ is of singular cardinality, then
$G \in K_p$ if and only if every pure subgroup $H$ of $G$ of smaller
cardinality than $G$ satisfies $H \in K_p$;
\item $G \not\in K_p$ if and only if $r_p^e(G) > 0$ whenever we add $|G|$ Cohen reals to the universe;
\item If $\rk(G) \geq \aleph_0$, then adding $|G|$ Cohen reals to
the universe adds a new member to $\Ext_p(G,\Z)$ preserving the old
ones.
\end{enumerate}
\end{lemma}

\proof
Let $G$ and $p$ be as stated. Part (i) is an easy application of the first author's Singular Compactness Theorem from \cite{Sh0}.\\
One implication of (ii) is trivial, hence assume that $G \not\in K_p$. By Lemma \ref{freerank} (ii) we may assume that $G$
does not have any pure subgroup of smaller rank than $G$ satisfying (ii). It is easily seen that the
rank $\delta=\rk(G)$ of $G$ must be uncountable. Thus $\delta>\aleph_0$ must be regular by (i). Let
$G=\bigcup\limits_{\alpha < \delta}G_{\alpha}$ be a filtration of $G$ by pure subgroups $G_{\alpha}$ of $G$ ($\alpha < \delta$). The claim now follows as in \cite{MeSh} repeating \cite[Theorem 9 and Theorem 10]{MeSh} (compare also Lemma \ref{meklerold}). The only difference is that in our situation the group $G$ is not almost free, hence we require in \cite[Theorem 10]{MeSh} that for every $\alpha$ in the stationary set $E$ there exists an element $a \in G_{\alpha+1}\backslash(G_{\alpha} + p^{\omega}G_{\alpha+1})$ which belongs to the $p$-adic closure of $G_{\alpha}+p^{\omega}G_{\alpha+1}$. This makes only a minor change in the proof of \cite[Theorem 10]{MeSh}.\\
Finally, (iii) follows similar to (ii) from the proof of \cite[Theorem 11]{MeSh}. The proof is therefore left to the reader.
\qed\\

Finally, we consider the {\it $p$-closure} of a pure subgroup $H$ of some
torsion-free abelian group $G$ which shall be needed in the proof of Theorem \ref{main1}.

\begin{definition}
Let $G$ be torsion-free and $H$ a pure subgroup of $G$. For every
prime $p \in \Pi$ the set
\[ {\cl}_{p}(G,H)=\{ x \in G: \textit{ for all } n \in \N \textit{ there
is } y_n \in H \textit{ such that } x -y_n \in p^nG \} \] is called
the {\rm $p$-closure of $H$}.
\end{definition}

We have a first easy lemma.

\begin{lemma}
\label{pclosure} Let $G$ be torsion-free and $H$ a pure subgroup of
$G$. Then the following hold for all primes $p \in \Pi$.
\begin{enumerate}
\item $H \subseteq \cl_p(G,H)$;
\item $\cl_p(G,H)$ is a pure subgroup of $G$;
\item $\cl_p(G,H)/H$ is $p$-divisible.
\end{enumerate}
\end{lemma}

\proof We fix a prime $p \in \Pi$. The first statement is trivial.
In order to prove (ii) assume that $mx \in \cl_p(G,H)$ for some $m
\in \N$ and $x \in G$. Then, for every $n \in \N$, there is $y_n \in
H$ such that $mx - y_n \in p^nG$, say $y_n=p^ng_n$ for some $g_n \in
G$. Without loss of generality we may assume that $(m,p)=1$. Hence
$1=km + lp^n$ for some $k,l \in \Z$. Thus
\[ x= kmx+lp^nx=kp^ng_n + ky_n + lp^nx \]
and hence $x-ky_n \in p^nG$ with $ky_n \in H$. Therefore $x \in
\cl_p(G,H)$ and (ii) holds. Finally, (iii) follows easily from (ii).
\qed\\

\section{Supercompact cardinals and large $p$-ranks}
In this section we shall assume that the existence of a supercompact
cardinal is consistent with $ZFC$. We shall then determine the
cardinal invariants $(r_0^ e(G), r_p^ e(G) : p \in \Pi)$ of
$\Ext(G,\Z)$ for every torsion-free abelian group in this model and
show that they are as large as possible. We start with a theorem
from \cite{MeSh} (see also \cite{Da}). Recall that for cardinals
$\mu, \gamma$ and $\delta$ we can define a partially ordered set
$Fn(\mu,\gamma,\delta)$ by putting
\[ Fn(\mu,\gamma,\delta)=\{ f : dom(f) \rightarrow \gamma : dom(f) \subseteq \mu, |dom(f)| < \delta \}. \]
The partrial order is given by $f \leq g$ if and only if $g \subseteq f$ as functions.

\begin{lemma}
Suppose $\kappa$ is a supercompact cardinal, $V$ is a model of $ZFC$
which satisfies $2^{\aleph_0}=\aleph_1$ and $\P=Fn(\mu,2,\aleph_1)
\times Fn(\rho,2,\aleph_0)$, where $\mu,\rho > \aleph_1$. Then $\P$
forces that every $\kappa$-free group is free.
\end{lemma}

\proof See \cite[Theorem 19]{MeSh}. \qed\\

As a corollary one obtains

\begin{lemma}
\label{model} If it is consistent with $ZFC$ that a supercompact
cardinal exists then both of the statements {\rm every
$2^{\aleph_0}$-free group is free and $2^{\aleph_0} < 2^{\aleph_1}$}
and {\rm every $2^{\aleph_0}$-free group is free and $2^{\aleph_0} =
2^{\aleph_1}$} are consistent with $ZFC$. Furthermore, if it is
consistent that there is a supercompact cardinal then it is
consistent that there is a cardinal $\kappa < 2^{\aleph_1}$ so that
if $\kappa$ Cohen reals are added to the universe then every
$\kappa$-free group is free.
\end{lemma}

\proof See \cite[Corollary 20]{MeSh}. \qed\\

We are now ready to prove our main theorem of this section working
in the model from Lemma \ref{model}. Assume that the existence of a
supercompact cardinal is consistent with $ZFC$. Let $V$ be any model
in which there exists a supercompact cardinal $\kappa$ such that the
weak diamond principle $\Diamond_{\lambda^+}^*$ holds for all
regular cardinals $\lambda \geq \kappa$. Now, we use Cohen forcing
to add $\kappa$ Cohen reals to $V$ to obtain a new model $\v$. Thus,
in $\v$ we still have $\Diamond_{\lambda^+}^*$ for all regular
$\lambda \geq \kappa$ and also $\Diamond_{\kappa}$ holds. Moreover,
we have $2^{\aleph_0}=\kappa$ and every $\kappa$-free group (of
arbitrary cardinality) is free by \cite{MeSh}.


\begin{theorem}
\label{main1} In any model $\v$ as described above, the following is
true for every non-free torsion-free abelian group $G$ and prime $p
\in \Pi$.

\begin{enumerate}
\item $r_0^e(G)=2^{\max\{\aleph_0,\fk_0(G)\}}$;
\item If $\fk_p(G)$ is finite, then $r_p^e(G)=\fk_p(G)$;
\item If $\fk_p(G)$ is infinite, then $r_p^e(G)=2^{\fk_p(G)}$.
\end{enumerate}
\end{theorem}

We would like to remark first that the above theorem shows that
$r^e_p(G)$ ($p \in \Pi_0$) is as large as possible for every
torsion-free abelian group $G$ in the model described above.
Moreover, Lemma \ref{existence} shows that every sequence of
cardinals $(\nu_p : p \in \Pi_0)$ not excluded by Theorem
\ref{main1} may be realized as the cardinal invariants of
$\Ext(H,\Z)$ for some torsion-free group $H$.

\proof Let $p \in \Pi_0$ be fixed. By Lemma \ref{freerank} we may
assume that $p^{\omega}G=0$ if $p \in \Pi$. Moreover, Lemma \ref{reduction} shows that also $\fk_p(G)=\rk(G)$ holds without loss of generality. We now prove the claim
by induction on the rank $\lambda=\rk(G)=\fk_p(G)$.\\

\underline{Case A:} \quad $\lambda$ is finite.\\
In this case we may assume without loss of generality that
$\Hom(G,\Z)=\{0\}$ since $G$ is of finite rank. If $p=0$, then
$r_0^e(G)=2^{\aleph_0}=\kappa$ follows from Lemma \ref{countabletf}
since $G$ is not free. Thus assume that $p > 0$. Then $r_p^e(G)$ is
the dimension of $\Hom(G,\Z/p\Z)$ as a vectorspace over $\Z/p\Z$.
Since $\Hom(G,\Z/p\Z)$ is the vectorspace dual of $G/pG$ it follows
that $r_p^e(G)=\rk(G)=\fk_p(G)$. Note that $G$ is $p$-reduced by assumption.\\

\underline{Case B:} \quad $\lambda=\aleph_0$.\\
Since $G$ is of countable rank it is well-known that there exists a
decomposition $G=G' \oplus F$ of $G$ where $F$ is a free group and
$G'$ satisfies $\Hom(G',\Z)=\{0\}$. Moreover, by assumption
$\fk_p(G)=\lambda=\aleph_0$, hence we obtain that
$\rk(G')=\fk_p(G')=\aleph_0$. Therefore,
$r_0^e(G)=r_0^e(G')=2^{\aleph_0}=\kappa$ follows from Lemma
\ref{countabletf}. If $p>0$ we conclude that $\Hom(G',\Z/p\Z)$ has
cardinality $2^{\aleph_0}$. Since $\Hom(G',\Z)=\{0\}$ it follows by
Lemma \ref{lemmaabschaetz} (ii) that
\[ 2^{\aleph_0} \leq r_p^e(G')=r_p^e(G) \leq 2^{\aleph_0} \]
and thus $r_p^e(G)=2^{\aleph_0}=\kappa$ for $p \in \Pi_0$.\\

\underline{Case C:} \quad $\aleph_0 < \lambda < \kappa$.\\
Note that $\kappa=2^{\aleph_0}=2^{\lambda}$, hence
$r_0^e(G)=2^{\aleph_0}=\kappa$ follows from Lemma
\ref{uncountabletf}. In fact, by induction hypothesis every
Whitehead group $H$ of size less than $\lambda$ has to be free
(because $r_p^e(H)=\fk_0(H)$) which suffices for Lemma
\ref{uncountabletf}. Now, assume that $ p \in \Pi$. By Lemma
\ref{lemmaabschaetz} (ii) we deduce
\[ r_p^e(G) \leq 2^{\fk_p(G)} \leq 2^{\lambda} = 2^{\aleph_0} \]
and hence it remains to prove that $r_p^e(G) \geq 2^{\aleph_0}$. The
proof is very similar to the proof of \cite[Theorem 8]{MeSh} (see
also Lemma \ref{meklerold}), hence we shall recall it only briefly.
Let $V$ be the ground model and $P$ be the Cohen forcing, i.e.
$P=P(\kappa \times \omega,2,\omega)=\{ h: \dom(h) \rightarrow
\{0,1\} : \dom(h) \textit{ is a finite subset of } \kappa \times
\omega \}$. If $\G$ is a $P$-generic filter over $V$ let
$\tilde{h}=\bigcup\limits_{g \in \G}g$. For notational reasons we
may also write $\v=V[\G]=V[\tilde{h}]$ for the extension model
determined by the generic filter $\G$. Let $A \subseteq
[\kappa]^{\leq \lambda}$ such that $G$ belongs to
$V[\tilde{h}\restriction_{A \times \omega}]$. Without loss of
generality we may assume that $\alpha \in A$ if and only if $\beta
\in A$ whenever $\alpha + \lambda=\beta + \lambda$. We shall prove
the claim by splitting the forcing. For each $\alpha$ such that
$\lambda\alpha \in \kappa\backslash A$ let
\[\tilde{f}_{\alpha} \in V[A \times \omega \cup [\lambda\alpha,
\lambda\alpha + \lambda) \times \omega]\] be a member of
$\Hom(G,\Z/p\Z)$ computed by
$\tilde{h}\restriction_{[\lambda\alpha,\lambda\alpha+\lambda)\times
\omega}$. Note that $\tilde{f}_{\alpha}$ exists by Proposition
\ref{forcingone} (iii). Then $\tilde{f}_{\alpha}$ is also computed
from
$\tilde{h}\restriction_{[\lambda\alpha,\lambda\alpha+\lambda)\times
\omega}$ over
$V[\tilde{h}\restriction_{\kappa\backslash[\lambda\alpha,\lambda\alpha+\lambda)
\times \omega}]$ and hence $\tilde{f}_{\alpha}$ is not equivalent to
any $f' \in
V[\tilde{h}\restriction_{\kappa\backslash[\lambda\alpha,\lambda\alpha+\lambda)\times
\omega}]$ modulo $\Hom(G,\Z)\varphi^p$. Thus the set of
homomorphisms $\{\tilde{f}_{\alpha} : \lambda\alpha < \kappa,
\lambda\alpha \not\in A \}$ exemplifies that $r_p^e(G) \geq
\kappa=2^{\aleph_0}$ and hence
\[2^{\lambda}=2^{\aleph_0}=\kappa \leq r_p^e(G) \leq 2^{\lambda}\] which shows
$r_p^e(G)=2^{\lambda}=2^{\fk_p(G)}$.\\

\underline{Case D:} \quad $\lambda \geq \kappa$.\\
Let $p \in \Pi_0$. We distinguish two subcases. Note that for $p \in \Pi$ our assumption
$\fk_p(G)=\rk(G)$ also implies that $\fk_0(G)=\rk(G)$ by Lemma \ref{freerank} (v).\\

\underline{Case D1:} \quad $\lambda \geq \kappa$ is regular.\\
The case $p=0$ follows as in \cite[Theorem XII 4.4]{EkMe}. Let
$\left< G_{\alpha} : \alpha < \lambda \right>$ be a filtration of
$G$ into pure subgroups $G_{\alpha}$ $(\alpha < \lambda)$ so that if
$G/G_{\alpha}$ is not $\lambda$-free, then $G_{\alpha+1}/G_{\alpha}$
is not free. Choose by \cite[Lemma 2.4]{EkHu} an associate free
resolution of $G$, i.e. a free resolution
\[ 0 \rightarrow K \overset{\Phi}{\rightarrow} F \rightarrow G \rightarrow 0 \]
of $G$ such that $F=\bigoplus\limits_{\alpha<\lambda}F_{\alpha}$ and
$K=\bigoplus\limits_{\alpha < \lambda}K_{\alpha}$ are free groups
such that $|F_{\alpha}| < \lambda$ and $|K_{\alpha}|<\lambda$ for
all $\alpha < \lambda$ and the induced sequences
\[ 0 \rightarrow \bigoplus\limits_{\beta < \alpha}K_{\beta}
\rightarrow \bigoplus\limits_{\beta < \alpha}F_{\beta} \rightarrow
G_{\alpha} \rightarrow 0\] are exact for every $\alpha < \lambda$.
Since $\fk_0(G)=\lambda$, the set $E=\{\alpha < \lambda :
G_{\alpha+1}/G_{\alpha} \text{ is not free } \}$ is stationary. For
any subset $E' \subseteq E$ let $K(E')=\bigoplus\limits_{\alpha \in
E'}K_{\alpha}$ and $G(E')=F/\Phi(K(E'))$. Then $\Gamma(G(E')) \geq
\tilde{E'}$ where $\Gamma(G(E'))$ is the $\Gamma$-invariant of
$G(E')$ (see \cite{EkMe} for details on the $\Gamma$-invariant).
Now, by assumption we have $\Diamond^*_{\lambda}$, hence we may
decompose $E$ into $\lambda$ disjoint stationary sets $E'_{\alpha}$,
each of which is non-small, i.e. $\Diamond_{\lambda}^*(E_{\alpha}')$
holds. Hence $G(E'_{\alpha})$ is not free since $\tilde{E'_{\alpha}}
\leq \Gamma(G(E'_{\alpha}))$ for every $\alpha < \lambda$. By
\cite{MeSh} we conclude that $G(E'_{\alpha})$ is not $\kappa$-free
and therefore has a non-free pure subgroup $H_{\alpha}$ of rank less
than $\kappa$. By induction hypothesis it follows that
$\Ext(H_{\alpha},\Z)\not= 0$ and hence also
$\Ext(G(E'_{\alpha}),\Z)\not= 0$. As in \cite[Lemma 1.1]{EkHu} (see
also \cite[Lemma XII 4.2]{EkMe}) there is an epimorphism
\[ \Ext(G,\Z) \rightarrow \prod\limits_{\alpha <
\lambda}\Ext(G(E_{\alpha}'),\Z) \rightarrow 0 \]
and it easily follows that $r_0^e(G) \geq 2^{\lambda}$ and hence
$r_0^e(G)=2^{\lambda}$ (compare \cite[Lemma 4.3]{EkMe}).\\

Now assume that $p>0$. Again, let $\left< G_{\alpha} : \alpha <
\lambda \right>$ be a filtration of $G$ into pure subgroups
$G_{\alpha}$ $(\alpha < \lambda)$ such that
$\Ext_p(G_{\alpha+1}/G_{\alpha},\Z)\not= 0$ if and only if
$\Ext_p(G_{\beta}/G_{\alpha},\Z)\not=0$ for some $\beta > \alpha$.
Fix $\alpha < \lambda$. We claim that $G/\cl_p(G,G_{\alpha})$ is not
free. By way of contradiction assume that $G/\cl_p(G,G_{\alpha})$ is
free. Hence $G=\cl_p(G,G_{\alpha}) \oplus F$ for some free group
$F$. Therefore, $G/G_{\alpha}=\left( \cl_p(G,G_{\alpha}) \oplus F
\right)/G_{\alpha}=\cl_p(G,G_{\alpha})/G_{\alpha} \oplus F$ is a
direct sum of a $p$-divisible group and a free group by Lemma
\ref{pclosure} (iii). It follows that $\fk_p(G) \leq \rk(G_{\alpha})
< \lambda$ contradicting the fact that $\fk_p(G)=\lambda$. By
\cite{MeSh} we conclude that $G/\cl_p(G,G_{\alpha})$ is not
$\kappa$-free since we are working in the model $\v$. Let
$G'/\cl_p(G,G_{\alpha}) \subseteq_* G/\cl_p(G,G_{\alpha})$ be a
non-free pure subgroup of $G/\cl_p(G,G_{\alpha})$ of size less than
$\kappa$. Then there exists $\alpha \leq \beta < \lambda$ such that
$G' \subseteq_* G_{\beta}$. By purity it follows that
$\cl_p(G,G_{\alpha}) \cap G_{\beta}=\cl_p(G_{\beta},G_{\alpha})$.
Hence
\[ G_{\beta}/{\cl}_{p}(G_{\beta},G_{\alpha}) = \left( G_{\beta} +
{\cl}_{p}(G,G_{\alpha}) \right) /{\cl}_{p}(G,G_{\alpha}) \] is
torsion-free but not free. Without loss of generality we may assume
that $\beta=\alpha+1$. Hence we may assume that for all $\alpha <
\lambda$ the quotient $G_{\alpha + 1}/\cl_p(G_{\alpha
+1},G_{\alpha})$ is a torsion-free non-free group. Note that
$G_{\alpha + 1}/\cl_p(G_{\alpha +1},G_{\alpha})$ is also $p$-reduced
since $\cl_p(G_{\alpha
+1},G_{\alpha})$ is the $p$-closure of $G_{\alpha}$ inside $G_{\alpha+1}$.\\
Since the cardinality of $G_{\alpha + 1}/\cl_p(G_{\alpha
+1},G_{\alpha})$ is less than $\lambda$ the induction hypothesis
applies. Hence $r_p^e(G_{\alpha + 1}/\cl_p(G_{\alpha
+1},G_{\alpha}))=2^{\fk_p(G_{\alpha + 1}/\cl_p(G_{\alpha
+1},G_{\alpha}))}$. We claim that $\Ext_p(G_{\alpha +
1}/\cl_p(G_{\alpha +1},G_{\alpha}),\Z) \not = \{0\}$ stationarily
often. If not, then there is a cub $C \subseteq \lambda$ such that
for all $\alpha < \beta \in C$ we have $\Ext_p(G_{\alpha +
1}/\cl_p(G_{\alpha +1},G_{\alpha}),\Z) = \{0\}$ and equivalently
$\Ext_p(G_{\beta}/G_{\alpha},\Z)=0$, hence
$\fk_p(G_{\beta})=\fk_p(G_{\alpha})$ by induction hypothesis. As in
\cite[Proposition XII 1.5]{EkMe} it follows that
$\Ext_p(G/G_{\alpha},\Z)=0$ for all $\alpha \in C$. This easily
contradicts the fact that $\fk_p(G)=\lambda$. It follows that
without loss of generality for every $\alpha < \lambda$ there exist
homomorphisms $h_{\alpha}^0, h_{\alpha}^1 \in
\Hom(G_{\alpha+1},\Z/p\Z)$ such that
\renewcommand{\labelenumi}{(\arabic{enumi})}
\begin{enumerate}
\item $h_{\alpha}^0\restriction_{\cl_p(G_{\alpha
+1},G_{\alpha})}=h_{\alpha}^1\restriction_{\cl_p(G_{\alpha
+1},G_{\alpha})}$;
\item There are no homomorphisms $g_{\alpha}^0,g_{\alpha}^1 \in
\Hom(G_{\alpha +1},\Z)$ such that
\begin{enumerate}
\item $g_{\alpha}^0\restriction_{\cl_p(G_{\alpha
+1},G_{\alpha})}=g_{\alpha}^1\restriction_{\cl_p(G_{\alpha
+1},G_{\alpha})}$ (or equivalently
$g_{\alpha}^0\restriction_{G_{\alpha}}=g_{\alpha}^1\restriction_{G_{\alpha}}$);
\item $h_{\alpha}^0=g_{\alpha}^0\varphi^p$ and
$h_{\alpha}^1=g_{\alpha}^1\varphi^p$. \end{enumerate}\end{enumerate}
To see this note that there is a homomorphism $\varphi:
G_{\alpha+1}/\cl_p(G_{\alpha+1},G_{\alpha}) \rightarrow \Z/p\Z$
which can not be factored by $\varphi_p$ since $\Ext_p(G_{\alpha +
1}/{\cl}_p(G_{\alpha +1},G_{\alpha}),\Z) \not = \{0\}$. Let
$h_{\alpha}^0=0$ and $h_{\alpha}^1$ be given by
\[ h_{\alpha}^1: G_{\alpha+1} \rightarrow
G_{\alpha+1}/{\cl}_p(G_{\alpha+1},G_{\alpha})
\overset{\varphi}\rightarrow \Z/p\Z.\] Then it is easy to check that
$h_{\alpha}^0$ and $h_{\alpha}^1$ are as required. In particular we
may assume that $h_{\alpha}^0=0$ for every $\alpha < \lambda$. An
immediate consequence is the following property (U).\\

\quad {   } \quad Let $f: G_{\alpha} \rightarrow \Z/p\Z$ and $g:
G_{\alpha} \rightarrow \Z$ such that $g\varphi_p=f$. Then there
exists\\
\quad (U) \quad $\tilde{f}: G_{\alpha+1} \rightarrow \Z/p\Z$
 such that $\tilde{f}\restriction_{G_{\alpha}}=f$ and
there is no homomorphism\\
{ \hspace*{0.5cm}} \quad $\tilde{g}:G_{\alpha+1} \rightarrow \Z$
satisfying both $\tilde{g}\restriction_{G_{\alpha}}=g$ and
$\tilde{g}\varphi_p=\tilde{f}$.\\

To see this, let $\hat{f}:G_{\alpha+1} \rightarrow \Z/p\Z$ be any
extension of $f$ which exists by the pure injectivity of $\Z/p\Z$.
If $\hat{f}$ is as required let $\tilde{f}=\hat{f}$. Otherwise let
$\hat{g}:G_{\alpha+1} \rightarrow \Z$ be such that
$\hat{g}\restriction_{G_{\alpha}}=g$ and $\hat{g}\varphi_p=\hat{f}$.
Put $\tilde{f}=\hat{f}-h_{\alpha}^1$. Then
$\tilde{f}\restriction_{G_{\alpha}}=f$. Assume that there exists
$\tilde{g}:G_{\alpha+1} \rightarrow \Z$ such that
$\tilde{g}\restriction_{G_{\alpha}}=g$ and
$\tilde{g}\varphi_p=\tilde{f}$. Choosing
$g_{\alpha}^1=\tilde{g}-\hat{g}$ we conclude
$-g_{\alpha}^1\varphi_p=h_{\alpha}^1$ contradicting (2). Note that
$h_{\alpha}^0=0$.\\
We now proceed exactly as in \cite[Proposition 1]{HHS} to show that
$r_p^e(G)=2^{\fk_p(G)}=2^{\lambda}$. We therefore recall the proof
only briefly and for simplicity we even shall assume that
$\Diamond_{\lambda}$ holds. It is an easy exercise (and therefore
left to the reader) to prove the result assuming the weak diamond
principle only. Assume that $r_p^e(G) =\sigma < 2^{\lambda}$ and let
$L=\{ f^{\alpha} : \alpha < \sigma \}$ be a complete list of
representatives of elements in $\Hom(G,\Z/p\Z)/\Hom(G,\Z)\varphi_p$.
Without loss of generality let $\{g_{\alpha}: G_{\alpha} \rightarrow
\Z: \alpha < \lambda \}$ be the Jensen functions given by
$\Diamond_{\lambda}$, hence for every homomorphism $g: G \rightarrow
\Z$ there exists $\alpha$ such that
$g\restriction_{G_{\alpha}}=g_{\alpha}$. We now define a sequence of
homomorphisms $\{f_{\alpha}^* : G_{\alpha} \rightarrow \Z/p\Z :
\alpha < \lambda\}$ such that the following hold.
\renewcommand{\labelenumi}{(\arabic{enumi})}
\begin{enumerate}
\item $f_{0}^*=f^0$;
\item $f^*_{\alpha}\restriction_{G_{\beta}}=f^*_{\beta}$ for all
$\beta < \alpha$;
\item If $f^*=\bigcup_{\alpha < \lambda}f^*_{\alpha}$, then
$f^*-f^{\alpha}$ is an element of $\Hom(G,\Z/p\Z)$ but not of
$\Hom(G,\Z)\varphi_p$.
\end{enumerate}
Suppose that $f^*_{\beta}$ has been defined for all $\beta <
\alpha$. If $\alpha$ is a limit ordinal, then we let
$f^*_{\alpha}=\bigcup_{\beta < \alpha}f^*_{\beta}$ which is a
well-defined homomorphism by (2). If $\alpha=\beta +1$ is a
successor ordinal, then we distinguish two cases. If
$f^*_{\beta}-f^{\beta}\restriction_{G_{\beta}}
\not=g_{\beta}\varphi_p$, let $f^*_{\alpha}:G_{\alpha} \rightarrow
\Z/p\Z$ be any extension of $f^*_{\beta}$ which exists since
$\Z/p\Z$ is pure injective and $G_{\beta} \subseteq_* G_{\alpha}$.
If $f^*_{\beta}-f^{\beta}\restriction_{G_{\beta}}
=g_{\beta}\varphi_p$, then (U) shows that there is a homomorphism
$\tilde{f}:G_{\alpha} \rightarrow \Z/p\Z$ extending
$f^*_{\beta}-f^{\beta}\restriction_{G_{\beta}}$ such that there is
no $\tilde{g}:G_{\beta+1} \rightarrow \Z$ with both extending
$g_{\beta}$ and $\tilde{g}\varphi_p=\tilde{f}$. Finally, put
$f_{\alpha}^*=\tilde{f} + f^{\alpha}\restriction_{G_{\alpha}}$ and
$f^*=\bigcup\limits_{\alpha < \lambda}f_{\alpha}^*$. It is now
straightforward to see that $f^*$ satisfies (3) and hence $f^*$
contradicts
the maximality of the list $L$.\\

\underline{Case D2:} \quad $\lambda \geq \kappa$ is singular.\\
First note that $\fk_p(G) > \kappa$ since $\kappa=2^{\aleph_0}$ is
regular. By induction on $\alpha < \lambda$ we choose subgroups
$K_{\alpha}$ of $G$ such that the following hold.
\renewcommand{\labelenumi}{(\arabic{enumi})}
\begin{enumerate}
\item $K_{\alpha}$ is a pure non-free subgroup of $G$;
\item $|K_{\alpha}| < \kappa$;
\item $K_{\alpha} \cap \sum\limits_{\beta < \alpha}K_{\beta}=\{ 0
\}$;
\item $\sum\limits_{\beta < \alpha}K_{\beta}$ is a pure subgroup
of $G$.
\end{enumerate}
Assume that we have succeeded in constructing the groups
$K_{\alpha}$ ($\alpha < \lambda$). Then \[ K=\sum\limits_{\beta <
\lambda}K_{\beta}=\bigoplus\limits_{\beta < \lambda}K_{\beta} \] is
a pure subgroup of $G$ and hence $r_p^e(G) \geq r_p^e(K)$ by Lemma
\ref{freerank} (ii). If $\Ext_p(K_{\alpha},\Z)=0$, then
$K_{\alpha}\in K_p$ follows by induction. Since $G$ is $p$-reduced
we obtain $K_{\alpha}\in K_0$ contradicting (1). Thus,
$\Ext_p(K_{\alpha},\Z) \not= \{0\}$ for every $\alpha < \lambda$
which implies that $r_p^e(K) \geq 2^{\lambda}$ since $\Ext(K,\Z)
\cong \prod\limits_{\alpha < \lambda}\Ext(K_{\alpha},\Z)$. It
therefore suffices to complete the construction of the groups
$K_{\alpha}$ $(\alpha < \lambda)$. Assume that $K_{\beta}$ for
$\beta < \alpha$ has been constructed. Let $\mu=(\kappa +
|\alpha|)^{< \kappa}$ which is a cardinal less than $\lambda$. Let
$H_{\alpha}$ be such that
\renewcommand{\labelenumi}{(\roman{enumi})}
\begin{enumerate}
\item $H_{\alpha} \subseteq_* G$;
\item $\sum\limits_{\beta < \alpha}K_{\beta} \subseteq H_{\alpha}$;
\item $|H_{\alpha}|=\mu$;
\item If $K \subseteq G$ is of cardinality less than $\kappa$, then
there is a subgroup $K' \subseteq_* G$ such that $H_{\alpha} \cap K
\subseteq K'$ and $K$ and $K'$ are isomorphic over $K \cap
H_{\alpha}$, i.e. there exists an isomorphism $\psi: K \rightarrow
K'$ which is the identity if restricted to $K \cap H_{\alpha}$.
\end{enumerate}
It is easy to see that $H_{\alpha}$ exists. Now, $G/H_{\alpha}$ is a
non-free group since $\fk_0(G)=\lambda$ and $p^{\omega}\left(
G/H_{\alpha} \right)=\{ 0 \}$. Hence \cite{MeSh} implies that there
is $K'_{\alpha} \subseteq G$ such that $\left( K_{\alpha}' +
H_{\alpha} \right) / H_{\alpha}$ is not free and $|K'_{\alpha}| <
\kappa$. Let $K_{\alpha}^0 \subseteq_* H_{\alpha}$ be as in (5),
i.e. $K_{\alpha}' \cap H_{\alpha} \subseteq K_{\alpha}^0$ and there
is an isomorphism $\psi_{\alpha}: K_{\alpha}' \rightarrow
K_{\alpha}^0$ which is the identity on $K_{\alpha}' \cap
H_{\alpha}$. Let $K_{\alpha}=\{ x - x\psi_{\alpha}
: x \in K_{\alpha}' \}$. Then $K_{\alpha}$ is as required. For instance
\[ K_{\alpha} \cong K_{\alpha}^0/\left( K_{\alpha}' \cap H_{\alpha} \right)
\cong K_{\alpha}'/\left( K_{\alpha}' \cap H_{\alpha} \right) \]
shows that $K_{\alpha}$ is not free. \qed\\


\begin{corollary}
In the model $\v$ let $\left< \mu_p : p \in \Pi_0 \right>$ be a
sequence of cardinals. Then there exists a torsion-free non-free
abelian group $G$ such that $r_p^e(G)=\mu_p$ for all $p \in \Pi_0$
if and only if
\begin{enumerate}
\item $\mu_0=2^{\lambda_0}$ for some infinite cardinal $\lambda_0$;
\item $\mu_p \leq \mu_0$ for all $p \in \Pi$;
\item $\mu_p$ is either
finite or of the form $2^{\lambda_p}$ for some infinite cardinal
$\lambda_p$. \end{enumerate}
\end{corollary}

\proof The proof follows easily from Lemma \ref{existence} and
Theorem \ref{main1}. \qed\\

\begin{corollary}
\label{maincor} In the model $\v$ let $\left< \mu_p : p \in \Pi_0
\right>$ be a sequence of cardinals. Then there exists a non-free
$\aleph_1$-free abelian group $G$ such that $r_p^e(G)=\mu_p$ for all
$p \in \Pi_0$ if and only if $\mu_p \leq \mu_0$ and
$\mu_p=2^{\lambda_p}$ for some infinite cardinal $\lambda_p$ for
every $p \in \Pi_0$.
\end{corollary}

\proof By Theorem \ref{main1} we only have to prove the existence
claim of the corollary. It suffices to construct $\aleph_1$-free
groups $G_p$ for $p \in \Pi_0$ such that
$r_p^e(G)=r_0^e(G)=2^{\aleph_0}=\kappa$ and $r_q^e(G)=0$ for all $p
\not=q \in \Pi$. Then $B=G_0^{(\lambda_0)} \oplus \bigoplus\limits_{p
\in \Pi} G_p^{(\lambda_p)}$ will be as required (see for instance the
proof of \cite[Theorem 3(b)]{HHS}). Fix $p \in \Pi_0$. From
\cite[Theorem XII 4.10]{EkMe} or \cite{SS2} it follows that there exists an $\aleph_1$-free non-free group $G_p$ of
size $2^{\aleph_0}$ such that $r_p^e(G_p)=2^{\aleph_0}=\kappa$ if $p
\in \Pi$. In \cite[Theorem XII 4.10]{EkMe} it is then assumed that
$2^{\aleph_0}=\aleph_1$ to show that also $r_0^e(G_p)=\kappa$.
However, since we work in the model $\v$ and $G_p$ is not free it
follows from Theorem \ref{main1} that $\fk_0(G_p) \geq \aleph_1$ and
hence $r_p^0(G_p)=2^{\fk_0(G_p)}=\kappa$. \qed\\

Recall that a reduced torsion-free group $G$ is called {\it
coseparable} if $\Ext(G,\Z)$ is torsion-free. By \cite{MeSh} it is
consistent that all coseparable groups are free. However, by
\cite[Theorem XII 4.10]{EkMe} there exist non-free coseparable
groups assuming $2^{\aleph_0}=\aleph_1$. Note that the groups
constructed in Lemma \ref{existence} are not reduced, hence do not
provide examples of coseparable groups.

\begin{corollary}
In the model $\v$ there exist non-free coseparable groups.
\end{corollary}

\proof Follows from Corollary \ref{maincor} letting
$\mu_0=2^{\aleph_0}$ and $\mu_p=0$ for all $p \in \Pi$. \qed


\section{A model close to $ZFC$}
In this section we shall construct a coseparable group which is not
free in a model of $ZFC$ which is very close to $ZFC$. As mentioned
in the previous Section $4$ it is undecidable in $ZFC$ if all
coseparable groups are free.\\
Let $\aleph_0 < \lambda$ be a regular cardinal and $S$ a stationary
subset of $\lambda$ consisting of limit ordinals of cofinality
$\omega$. We recall the definition of a ladder system on $S$ (see
for instance \cite[page 405]{EkMe}).

\begin{definition}
{\rm A} ladder system $\bar{\eta}$ on $S$ {\rm is a family of
functions $\bar{\eta}=\left< \eta_{\delta} : \delta \in S \right>$
such that $\eta_{\delta}: \omega \rightarrow \delta$ is strictly
increasing with $\sup(\rg(\eta_{\delta}))=\delta$, where
$\rg(\eta_{\delta})$ denotes the range of $\eta_{\delta}$. We call
the ladder system} tree-like {\rm if for all $\delta, \nu \in S$ and
every $n, n\in \omega$, $\eta_{\delta}(n)=\eta_{\nu}(m)$ implies
$n=m$ and $\eta_{\delta}(k)=\eta_{\nu}(k)$ for all $k \leq n$.}
\end{definition}

In order to construct almost-free groups one method is to use
$\kappa$-free ladder systems.
\begin{definition}
Let $\kappa$ be an uncountable regular cardinal. The ladder
 system $\bar{\eta}$ is called {\rm $\kappa$-free} if for every
 subset $X \subseteq S$ of cardinality less than $\kappa$ there is
 a sequence of natural numbers $\left< n_{\delta} : \delta \in X \right>$
 such that
 \[ \left< \{\eta_{\delta}\restriction_l : n_{\delta} < l < \omega\} :
 \delta \in X \right> \]
 is a sequence of pairwise disjoint sets. \end{definition}

Finally, recall that a stationary set $S \subseteq \lambda$ with $\lambda$
uncountable regular is called {\it non-reflecting}
if $S \cap \kappa$ is not stationary in $\kappa$ for every $\kappa < \lambda$
with $\cf(\kappa) > \aleph_0$.

\begin{theorem}\label{main2}
Let $\mu$ be an uncountable strong limit cardinal such that
$\cf(\mu)=\omega$ and $2^{\mu}=\mu^+$. Put $\lambda=\mu^+$ and
assume that there exists a $\lambda$-free tree-like ladder system on
a non-reflecting stationary subset $S \subseteq \lambda$. If
$\Pi=\Pi_0 \cup \Pi_1$ is a partition of $\Pi$ into disjoint subsets
$\Pi_0$ and $\Pi_1$, then there exists an almost-free group $G$ of
size $\lambda$ such that
\begin{enumerate}
\item $r_0^e(G)=2^{\lambda}$;
\item $r_p^e(G)=2^{\lambda}$ if $p \in \Pi_0$;
\item $r_p^e(G)=0$ if $p \in \Pi_1$.
\end{enumerate}
\end{theorem}

\proof Let $\bar{\eta}=\left< \eta_{\delta} : \delta \in S \right>$
be the $\lambda$-free ladder system where $S$ is a stationary non-reflecting
subset
of $\lambda$ consisting of ordinals less than $\lambda$ of
cofinality $\omega$. Without loss of generality we may assume that
$S=\lambda$. Let $\pr : \mu^2 \rightarrow \mu$ be a pairing
function, hence $\pr$ is bijective and if $\alpha \in \mu$ then we
shall denote by $(\pr_1(\alpha),\pr_2(\alpha))$ the unique pair
$(\beta,\gamma) \in \mu^2$ such that
$\pr(\pr_1(\alpha),\pr_2(\alpha))=\alpha$. Let $L$ be the free
abelian group
\[ L=\bigoplus\limits_{\alpha < \mu} \Z x_{\alpha} \]
generated by the independent elements $x_{\alpha}$ ($\alpha < \mu$).
For notational simplicity we may assume that $\Pi_1 \not= \emptyset$
and let $\left< (p_{\beta},f_{\beta}) : \beta < \lambda \right>$ be
a listing of all pairs $(p,f)$ with $p \in \Pi_1$ and $f \in
\Hom(L,\Z/p\Z)$. Recall that $\lambda=2^{\mu}$. By induction on
$\beta < \lambda$ we shall choose triples
$(g_{\beta},\nu_{\beta},\rho_{\beta})$ such that the following
conditions hold.
\renewcommand{\labelenumi}{(\arabic{enumi})}
\begin{enumerate}
\item $g_{\beta} \in \Hom(L,\Z)$;
\item $f_{\beta}=g_{\beta}\varphi_p$ where $\varphi_p: \Hom(L,\Z) \rightarrow \Hom(L,\Z/p\Z)$ is the canonical map;
\item $\nu_{\beta},\rho_{\beta} : \omega \rightarrow \mu$ such that
$\eta_{\beta}(n)=\pr_1(\nu_{\beta}(n))=\pr_1(\rho_{\beta}(n))$;
\item For all $\delta \leq \beta$ there exists $n=n(\delta, \beta) \in \omega$ such that
for all $m \geq n$ we have
$g_{\delta}(x_{\nu_{\beta}(m)})=g_{\delta}(x_{\rho_{\beta}(m)})$;
\item For all $\delta < \beta$ there exists $n=n(\delta, \beta) \in \omega$ such
that
for some sequence $\left< b^{\delta,\beta}_m : m \in [n,\omega) \right>$
of natural numbers we have $\left(\prod\limits_{p \in \Pi_1 \cap m}p
\right)b^{\delta,\beta}_{m+1}=b^{\delta,\beta}_m +
g_{\beta}(x_{\nu_{\delta}(m)}) - g_{\beta}(x_{\rho_{\delta}(m)})$
for all $m \geq n$;
\item $\nu_{\beta}(m) \not= \rho_{\beta}(m)$ for all $m \in \omega$.
\end{enumerate}
Fix $\beta < \lambda$ and assume that we have constructed
$(g_{\delta},\nu_{\delta},\rho_{\delta})$ for all $\delta < \beta$.
Choose a function $h_{\beta}: \beta \rightarrow \omega$ such that
$h_{\beta}(\delta) > p_{\delta}$ for all $\delta < \beta$ and
\begin{equation}
\label{disjoint} \left< \{ \eta_{\delta}\restriction_l : l \in
[h_{\beta}(\delta),\omega)\} : \delta < \beta \right> \end{equation}
is a sequence of pairwise disjoint sets. Note that such a choice is
possible since the ladder system $\bar{\eta}$ is $\lambda$-free by
assumption. Moreover, by (3) the pairing function $\pr$ implies that
also \begin{equation} \label{almostfree} \left< \{
\nu_{\delta}\restriction_l, \rho\restriction_l : l \in
[h_{\beta}(\delta),\omega)\} : \delta < \beta \right> \end{equation}
is a sequence of pairwise disjoint sets. Now, we choose the function $g_{\beta}$ such
that (2) and (5) hold.
For $\delta < \beta$ let $n=n(\delta,\beta)=h_{\beta}(\delta)$. Since $L$
is free we may choose first $g_{\beta}(x_{\alpha})$ satisfying
$g_{\beta}(x_{\alpha}) + p_{\beta}\Z = f_{\beta}(x_{\alpha})$ for
every $\alpha$ such that $\pr_1(\alpha) \not= \eta_{\delta}(l)$ for
all $\delta < \beta$ and $l \geq n(\delta,\beta)$, that is to say for
those $\alpha$ such that $x_{\alpha}$ does not appear in (5).
Secondly, for $\delta < \beta$, we choose by induction on $m \geq
n(\delta,\beta)$ integers $b^{\delta,\beta}_{m+1}$ such that
\[ 0 + p_{\beta}\Z = b^{\delta,\beta}_{m+1} +
f_{\beta}(x_{\nu_{\delta}(m)}) - f_{\beta}(x_{\rho_{\delta}(m)}) +
p_{\beta}\Z \] and then choose $g_{\beta}(x_{\nu_{\delta}(m)})$ and
$g_{\beta}(x_{\rho_{\delta}(m)})$ such that (5) holds for $\delta$.
Note that this inductive process is possible by the choice of
$h_{\beta}$ and condition (\ref{disjoint}).\\
Finally, let $\beta=\bigcup\limits_{n \in \omega}A_n$ be the union
of an increasing chain of sets $A_n$ such that $|A_n|< \mu$ (recall
that we have assumed without loss of generality that $S=\lambda$, so
$\beta$ is of cofinality $\omega$). By induction on $n < \omega$ we
now may choose $\rho_{\beta}(n)$ and $\nu_{\beta}(n)$ as distinct
ordinals such that \begin{itemize}
\item $\rho_{\beta}(n),\nu_{\beta}(n) \in \mu$
\item $\rho_{\beta}(n),\nu_{\beta}(n) \not\in \{ \nu_{\beta}(m),
\rho_{\beta}(m) : m < n \}$
\item
$\pr_1(\rho_{\beta}(n))=\pr_1(\nu_{\beta}(n))=\eta_{\beta}(n)$;
\item $\left< g_{\delta}(x_{\nu_{\beta}(n)}) : \delta \in A_n
\right> = \left< g_{\delta}(x_{\rho_{\beta}(n)}) : \delta \in A_n
\right>$.
\end{itemize}
Hence (3), (4) and (6) hold and we have carried on the induction.
Now, let $G$ be freely generated by $L$ and $\{ y_{\beta,n} : \beta
< \lambda, n \in \omega \}$ subject to the following relations for
$\beta < \lambda$ and $n \in \omega$.
\[ \left( \prod\limits_{p \in \Pi_1 \cap n} p \right) y_{\beta,n+1}
= y_{\beta,n}+x_{\nu_{\beta}(n)} - x_{\rho_{\beta}(n)} .\] Then $G$
is a torsion-free abelian group of size $\lambda$. Moreover, since
the ladder system $\bar{\eta}$ is $\lambda$-free and $S$ is
stationary but not reflecting it follows by standard calculations
using (\ref{almostfree}) that $G$ is almost-free but not free (see for instance \cite{EkSh2}).
It
remains to prove that (i), (ii) and (iii) of the Theorem hold. For
$\beta < \lambda$ let
\[ G_{\beta}= \left< L, y_{\delta,n} : \delta < \beta, n \in \omega
\right>_* \subseteq_* G
\] so that $G=\bigcup\limits_{\beta < \lambda}G_{\beta}$ is the
union of the continuous increasing sequence of pure subgroups
$G_{\beta}$ ($\beta < \lambda$). We start by proving (iii). Thus let
$p \in \Pi_1$ and choose $f \in \Hom(G,\Z/p\Z)$. By assumption there
is $\beta < \lambda$ such that
$(p,f\restriction_L)=(p_{\beta},f_{\beta})$. Inductively we shall
define an increasing sequence of homomorphisms $g_{\beta,\gamma}:
G_{\gamma} \rightarrow \Z$ for $\gamma \geq \beta$ such that
$g_{\beta,\gamma}\varphi_p=f\restriction_{G_{\gamma}}$. For
$\gamma=\beta$ we choose $n(\delta,\beta)$ and $\left< b^{\delta,\beta}_m : m
\in [n(\delta,\beta),\omega) \right>$ as in (5) for $\delta < \beta$. We let
$g_{\beta,\beta}\restriction_L=g_{\beta}$ where $g_{\beta}$ is
chosen as in (1). Moreover, put
$g_{\beta,\beta}(y_{\delta,m})=b_m^{\delta,\beta}$ for $m \in
[n(\delta,\beta),\omega)$ and $\delta < \beta$. By downwards induction we
chose $g_{\beta,\beta}(y_{\delta,m})$ for $m < n(\delta\beta)$, $\delta <
\beta$. It is easily seen that $g_{\beta,\beta}$ is as required,
i.e. satisfies $g_{\beta,\beta}\varphi_p=f\restriction_{G_{\beta}}$.
Now, assume that $\gamma > \beta$. If $\gamma$ is a limit ordinal,
then let $g_{\beta,\gamma}=\bigcup\limits_{\beta \leq \epsilon <
\gamma}g_{\beta,\epsilon}$. If $\gamma=\epsilon+1$, then (4) implies
that there is $n(\beta,\epsilon) < \omega$ such that
$g_{\beta}(x_{\nu_{\epsilon}(m)})=g_{\beta}(x_{\rho_{\epsilon}(m)})$ for
all $m \in [n(\beta,\epsilon),\omega)$. Therefore, putting
$g_{\beta,\gamma}\restriction_{G_{\epsilon}}=g_{\beta,\epsilon}$ and
$g_{\beta,\gamma}(y_{\epsilon,m})=0$ for $m \in [n(\beta,\epsilon),\omega)$
and determing $g_{\beta,\gamma}(y_{\epsilon,m})$ by downward
induction for $m < n(\beta,\epsilon)$ we obtain $g_{\beta,\gamma}$ as
required. Finally, let $g=\bigcup\limits_{\gamma \geq
\beta}g_{\beta,\gamma}$ which satisfies $g\varphi_p=f$. Since $f$
was chosen arbitrary it follows that
$\Hom(G,\Z/p\Z)=\Hom(G,\Z)\varphi_p$ for
all $p \in \Pi_1$ and hence $r_p^e(G)=0$ for $p \in \Pi_1$.\\
We now turn to $p \in \Pi_0$. By definition of $G$ it follows that
every homomorphism $\psi:L \rightarrow \Z$ has at most one extension
to a homomorphism $\psi':G \rightarrow \Z$. Thus $|\Hom(G,\Z)|\leq
2^{\mu}$. However, for every $\beta < \lambda$, any homomorphism
$\psi: G_{\beta} \rightarrow \Z/p\Z$ has more than one extension to
a homomorphism $\psi':G_{\beta+1} \rightarrow \Z/p\Z$ and hence
$|\Hom(G,\Z/p\Z)|=2^{\lambda} > 2^{\mu}$. Consequently,
$r_p^e(G)=2^{\lambda}$. Similarly, it follows that $r_0^e(G)=2^{\lambda}$ which
finishes the proof. \qed\\


\begin{corollary}
Let $\mu$ be an uncountable strong limit cardinal such that
$\cf(\mu)=\omega$ and $2^{\mu}=\mu^+$. Put $\lambda=\mu^+$ and
assume that there exists a $\lambda$-free ladder system on a
stationary subset $S \subseteq \lambda$. Then there exists an
almost-free non-free coseparable group of size $\lambda$.
\end{corollary}

\proof Follows from Theorem \ref{main2} letting $\Pi_0=\emptyset$
and $\Pi_1=\Pi$. \qed

\goodbreak

\bigskip

\end{document}